\documentclass[11pt,letterpaper]{article}
\usepackage[pass]{geometry}

\usepackage{palatino}
\usepackage{mathpazo}
\usepackage{stmaryrd}

\usepackage{hyperref}
\hypersetup{pdfpagemode=UseNone}

\usepackage{amsfonts}   
\usepackage{amssymb}  
\usepackage{amsmath}
\usepackage{latexsym}
\usepackage{amsthm}
\usepackage{sectsty}

\usepackage{color}
\usepackage{graphicx}    
\usepackage{mathrsfs}
\usepackage{bbm,bm}

\theoremstyle{definition}

\def\seq{\operatorname{\hbox{\sc Seq}}}
\def\set{\operatorname{\hbox{\sc Set}}}
\def\cyc{\operatorname{\hbox{\sc Cyc}}}

\begin{document}

\title{Combinatorial interpretation and proof \\of Glaisher--Crofton identity}

\author{P. Blasiak\thanks{Corresponding author: \href{mailto:pawel.blasiak@ifj.edu.pl}{pawel.blasiak@ifj.edu.pl}} , G. H. E. Duchamp, A. Horzela \& K. A. Penson}

\date{}

\maketitle

\begin{abstract}
We give a purely combinatorial proof of the Glaisher-Crofton identity which derives from the analysis of discrete structures generated by iterated second derivative. The argument illustrates utility of symbolic and generating function methodology of modern enumerative combinatorics and their applications to computational problems. 
\end{abstract}

\vspace{0.1cm}


%
%
%

\section{Introduction}

Many computational problems involve action of complex expressions in derivatives on functions. A typical example is exponential of the hamiltonian acting on the initial condition
which is a formal solution to evolution equation. Application of the latter range from classical heat and diffusion theory, financial mathematics and economy to quantum field theory. Hence practical interest in operational formulas enabling explicit evaluation of such expressions. Methods used to this effect usually involve operator and special function techniques, integral transforms, umbral calculus methods, etc. (see comprehensive review of the subject~\cite{DaOtToVa97,RiTa09}). In this paper we develop another approach based on modern combinatorial methods of analysis and enumeration of structures via generating functions~\cite{FlSe09,BeLaLe98,Wi06}.

The most known operational identities that involve exponential of derivative are formulas for the shift and dilation operators
\begin{eqnarray}\label{shift-dilation}
\exp\left(\lambda\,\tfrac{d}{dx}\right)F(x)=F(x+\lambda)&\text{and}&\exp\left(\lambda\,x\tfrac{d}{dx}\right)F(x)=F(e^\lambda x)\,,
\end{eqnarray}
where $F(x)$ is an arbitrary function\footnote{Here, we leave subtle problems of convergence aside and consider $F(x)$ as a formal power series in one variable $x$.}.
They are a special case of the general closed-form operational expression
\begin{eqnarray}\label{q-v}
\exp\left[\lambda\left(q(x)\,\tfrac{d}{dx} + v(x)\right)\right]F(x) = g(\lambda, x)\,
\cdot \,F(T(\lambda,x))\,,
\end{eqnarray}
where functions $T(\lambda,x)$ and $g(\lambda,x)$ are specified by the set of equations
\begin{eqnarray}
\frac{\displaystyle{\partial}T(\lambda,x)}{\displaystyle{\partial}\lambda} =  q(T(\lambda,x))\ , &&T(0,x) = x\ ,\\\frac{\displaystyle{\partial}g(\lambda,x)}{\displaystyle{\partial}\lambda} =  v(T(\lambda,x))\cdot g(\lambda,x)\ , &&g(0,x) = 1\ .
\end{eqnarray}
See \cite{DaOtToVa97,BlHoPeDuSo05,RiTa09} for the proof based on operator techniques, and~\cite[Sect.~6]{BlFl11} for a recently developed combinatorial approach.
 
Note that formulas Eqs.~(\ref{shift-dilation}) and~(\ref{q-v}) are valid for any function $F(x)$. However, this is a very unique situation which holds only for exponential of an expressions linear in the first derivative. For the second derivative closed-form of expressions are not known and the best what can be done is evaluation on specific functions. There are only a few explicit examples which include formula for the exponential generating function of Hermite polynomials~\cite{SrMa84}
\begin{eqnarray}\label{Hermite}
\exp\left(-\tfrac{1}{2}\tfrac{d^2}{d x^2}\right)\,\exp\left(xt\right)=\sum_{n\geqslant 0}He_n(x)\frac{t^n}{n!}=\exp\left(xt-\tfrac{1}{2}t^2\right)\,,
\end{eqnarray}
and the Glaisher-Crofton identity~\cite{Cr79,DaKhRi08,Da00}
\begin{eqnarray}\label{Glaisher}
\exp\left(\alpha\,\tfrac{d^2}{d x^2}\right)\,\exp\left(-x^2\right)=\frac{1}{\sqrt{1+4\alpha}}\exp\left(-\frac{x^2}{1+4\alpha}\right)\,.
\end{eqnarray}

Formulas of this type are usually derived using integral representations in the complex domain\footnote{For example, to derive Eq.~(\ref{Glaisher}) one may use integral representation $\exp\left(-x^2\right)=\int_{-\infty}^{+\infty}\exp\left(-\xi^2+2i\xi x\right)d\xi$, see~\cite{DaKhRi08,Da00}.}. However, in this paper we demonstrate that it is also possible to prove these identities on the basic algebraic level by analysing combinatorial structures generated by the action of second derivatives. In this approach functions are treated as generating functions enumerating simple combinatorial objects (like sets of subsets, cycles, sequences, etc.) and consequently expressions in derivatives transform objects in the class into richer structures whose generating functions can be quickly identified with the methods of symbolic combinatorics. Our goal in this paper is to develop and promote general combinatorial methodology for solving computational problems which, in many cases, provides better insight into algebraic and analytic manipulations (see e.g.~\cite{BlFl11}). We illustrate this novel approach by explaining combinatorial meaning of the exponential of second derivative and use this interpretation to derive Eqs.~(\ref{Hermite}) and (\ref{Glaisher}).

\section{Combinatorics of derivatives}

Action of derivative on a function can be understood as a transformation of combinatorial structure. One then requires that on the level of generating functions it corresponds to derivation. In the following we formalise this intuition by describing the relevant constructions and develop a broader picture which includes higher derivatives and their exponentials. This framework will be illustrated by a simple combinatorial proof and interpretation of Taylor's formula and the identity Eq.~(\ref{Hermite}).

\subsection{Generating functions, first derivative and Taylor's formula}

Let us consider a combinatorial class $\mathscr{F}$ which is defined as a denumerable collection of objects built of atoms represented by $\mathcal{X}$ according to some well specified procedure. A typical combinatorial problem consists in enumeration of objects in $\mathscr{F}$ according to the size $|\cdot|:\mathscr{F}\rightarrow\mathbb{N}$ which is usually the number of atoms. In other words, one seeks the sequence $f_n=\#\{\phi\in\mathscr{F}:|\phi|=n\}$ which counts the number of objects comprised of exactly $n$ atoms. This sequence can be encoded in a generating function
\begin{eqnarray}
F(x)=\sum_{\phi\in\mathscr{F}}x^{|\phi|}=\sum_{n\geqslant0}f_n\,x^n\,
\end{eqnarray}
which is a convenient tool for enumeration of complex structures via the so called transfer rules. The latter translate combinatorial constructions into algebraic manipulations of the corresponding generating functions (see~\cite{FlSe09,Wi06,BeLaLe98} for a comprehensive treatment of the subject and Appendix~\ref{Appendix} for a quick extract of a few transfer rules used in this paper). In the following we will briefly review combinatorics of the first derivative and recall a simple combinatorial interpretation of Taylor's formula.

Here we will be concerned with the derivative operation $yD_x$ acting on the class $\mathscr{F}$ which consists in: \emph{"selecting in all possible ways a single atom of type $\mathcal{X}$ in each element of $\mathscr{F}$ and replacing it with atom of a new type $\mathcal{Y}$"}. In other words, one may think of the new class $yD_x\mathscr{F}$ as formed of all structures taken from $\mathscr{F}$ in which one of the atoms $\mathcal{X}$ gets 'repainted' into a new colour $\mathcal{Y}$. Since each structure in $\mathscr{F}$ built of $n$ atoms of type $\mathcal{X}$ gives rise to $n$ new ones with $(n-1)$ atoms $\mathcal{X}$ and a single $\mathcal{Y}$, then the generating function enumerating objects in the new class $yD_x\mathscr{F}$ is given by
\begin{eqnarray}
G(x,y)=\sum_{n,k\geqslant0}g_{n,k}\,x^ny^k=\sum_{n\geqslant1}n\cdot f_n\ x^{n-1}y\,,
\end{eqnarray} 
where $g_{n,k}$ counts objects according to the number of $\mathcal{X}$'s and $\mathcal{Y}$'s respectively. This substantiates in the standard transfer rule
\begin{eqnarray}\label{D-transfer}
\text{{\textsl{Selecting singleton:}}}\qquad\mathscr{G}=yD_x\mathscr{F}&\Longrightarrow&G(x,y)=y\tfrac{d}{d x}F(x)\,.
\end{eqnarray}

Now, let us consider the $k$-th derivative $\tfrac{1}{k!}(yD_x)^k$ acting on $\mathscr{F}$. Combinatorially it means: \emph{"select in all possible ways an unordered collection of $k$ atoms of type $\mathcal{X}$ and replace (repaint) them by atoms of type $\mathcal{Y}$"}. Clearly, for each structure of size $n$ in $\mathscr{F}$ we have $\binom{n}{k}$ possible choices, and hence the generating function of the new class $\tfrac{1}{k!}(yD_x)^k\mathscr{F}$ evaluates to
\begin{eqnarray}
G(x,y)=\sum_{n,l\geqslant0}g_{n,l}\,x^ny^l=\sum_{n\geqslant1}\tiny{\binom{n}{k}}\, f_n\ \ x^{n-k}y^k\,.
\end{eqnarray} 
In consequence we get the following transfer rule
\begin{eqnarray}\label{Dk-transfer}
\text{{\textsl{Selecting $k$-subset:}}}\qquad\mathscr{G}=\tfrac{1}{k!}(yD_x)^k\mathscr{F}&\Longrightarrow&G(x,y)=\tfrac{1}{k!}y^k\tfrac{d^k}{d x^k}F(x)\,,
\end{eqnarray}
which gives combinatorial interpretation of $k$-th derivative elevated to the level of combinatorial structures. 

This brings an interesting perspective on derivative of a function which we will develop throughout the paper. Namely, one may think of a function $F(x)$ as the generating function of some combinatorial class $\mathscr{F}$. Then differentiation yields a new generating function which enumerates objects in the new class comprised of structures taken from $\mathscr{F}$ in which some of the atoms $\mathcal{X}$ were replaced with $\mathcal{Y}$s (how many are replaced depends on the order of derivative). Hence, derivative of a function $F(x)$ can be understood as a well defined combinatorial transformation of the associated combinatorial class $\mathscr{F}$ in a sense that on the level of generating functions it corresponds to simple differentiation (cf. Eqs.~(\ref{D-transfer}) and (\ref{Dk-transfer})).

For illustration of this this view point let us recall the usual Taylor's formula
\begin{eqnarray}\label{Taylor}
\sum_{k\geqslant0}\tfrac{y^n}{n!}\,F^{(n)}(x)=F(x+y)\,.
\end{eqnarray}
Surprisingly, it admits a transparent combinatorial interpretation (see~\cite[Note III.31]{FlSe09} or~\cite[Note 3]{BlFl11}). To see this one observes that the l.h.s. is the sum of derivatives $\exp \left(yD_x\right)\equiv\sum_{k\geqslant0}\tfrac{1}{k!}(yD_x)^k$ applied to function $F(x)$ which can be considered as the generating function of a class of objects built of atoms $\mathcal{X}$. Then from our previous discussion the exponential $\exp \left(yD_x\right)$ corresponds to: \emph{"selecting in all possible ways an \emph{arbitrary} number of $\mathcal{X}$'s and replacing them by $\mathcal{Y}$'s"} (since the sum contains all derivatives we may choose a subset of arbitrary cardinality). On the other hand, this is the same as substituting each atom $\mathcal{X}$ either with atom $\mathcal{X}$ (which makes nothing) or atom $\mathcal{Y}$ (which means replacement). Hence we have the following combinatorial equivalence
\begin{eqnarray}\label{Taylor-spec}
\exp \left(yD_x\right)\mathscr{F}=\mathscr{F}\circ(\mathcal{X}+\mathcal{Y})\,,
\end{eqnarray}
which on the level of generating functions, by virtue of Eq.~(\ref{Dk-transfer}) and the transfer rule for substitutions (see Appendix~\ref{Appendix}, Eq.~(\ref{substitution})), directly translates into Eq.~(\ref{Taylor}). Hence from the combinatorial point of view Taylor's formula is a simple manifestation of the following transfer rule (cf. Eq.~(\ref{Taylor-spec}))
\begin{eqnarray}\label{expD}
\text{{\textsl{Selecting arbitrary subset:}}}\qquad\mathscr{G}=\exp \left(yD_x\right)\mathscr{F}&\Longrightarrow&G(x,y)=F(x+y)\,,
\end{eqnarray}
which applies to any combinatorial class $\mathscr{F}$ and its generating function $F(x)$. This is a typical example of combinatorial methodology which draws on the fact that in many cases the same combinatorial structure allows different specifications.


\subsection{Second derivative and Hermite polynomials}

Combinatorial interpretation of the second derivative can be developed along similar lines. From the above we know that the 2-nd derivative $\tfrac{1}{2}(yD_x)^2$ acting on $\mathscr{F}$ consists in \emph{"selecting in all possible ways an unordered pair of $\mathcal{X}$ atoms and then replacing chosen $\mathcal{X}$'s by atoms of type $\mathcal{Y}$"}. We will call such (unordered) pair a \emph{doubleton}. Clearly, for an object composed of $n$ atoms this can be done in $\tfrac{n(n-1)}{2}$ ways, which agrees with the algebraic identity $\tfrac{1}{2}(yD_x)^2\ x^n=\tfrac{n(n-1)}{2}\,x^{n-2}y^2$.

More generally, by iterating $k$ times one picks up a sequence of doubletons in the original structure. Hence we define the following construction
\begin{eqnarray}\label{yD2k}
\text{{\textsl{Selecting $k$-subset of doubletons:}}}\qquad\mathscr{G}=\tfrac{1}{k!}\left(\tfrac{1}{2}\left(yD_x\right)^2\right)^k\mathscr{F}&\Longrightarrow&\tfrac{1}{2^kk!}\,y^{2k}\tfrac{d^{2k}}{dx^{2k}}F(x)\,,
\end{eqnarray}
which consists in: \emph{"selecting in all possible ways a set of $k$ unordered pairs (doubletons) of $\mathcal{X}$ atoms and replacing them by $\mathcal{Y}$ atoms"}. Note that, as in Eq.~(\ref{Dk-transfer}), we deem order in the sequence irrelevant by introducing factor $\tfrac{1}{k!}$ in front of the iterated derivative (hence the 'set' and not 'sequence' in the description). For a quick check of this specification we observe that the coefficient on the r.h.s. of the identity
\begin{eqnarray}
\tfrac{1}{k!}\left(\tfrac{1}{2}\left(yD_x\right)^2\right)^kx^n=\tfrac{1}{k!}\,\tfrac{n(n-1)}{2}\cdot\tfrac{(n-2)(n-3)}{2}\cdot ...\cdot \tfrac{(n-2k+2)(n-2k+1)}{2}\ x^{n-2k}y^{2k}
\end{eqnarray}
coincides with the number of possible choices of a set of $k$ unordered pairs from the set of $n$ objects, and hence by linearity we establish correctness of the description and transfer rule of Eq.~(\ref{yD2k}).

Now, we are in position to give interpretation of the exponential of second derivative, i.e.
\begin{eqnarray}\label{expyD2}
\text{{\textsl{Selecting arbitrary subset of doubletons:}}}\qquad\mathscr{G}=\exp \left(\tfrac{1}{2}(yD_x)^2\right)\mathscr{F}\equiv\sum_{k\geqslant0}\tfrac{1}{k!}(\tfrac{1}{2}yD_x)^{2k}\mathscr{F}\,,
\end{eqnarray}
whose combinatorial meaning comes down to: \emph{"selecting in all possible ways an \emph{arbitrary} subset of (unordered) pairs of $\mathcal{X}$ atoms and replacing each chosen $\mathcal{X}$ by atom of type $\mathcal{Y}$"}. This is the sum of derivative operations of the type Eq.~(\ref{yD2k}) which on the level of generating functions is the sum of derivatives. Unfortunately it does not close to any neat expression like Eqs.~(\ref{expD}) or (\ref{Taylor}); indeed, Taylor's formula does not generalise in a straightforward manner. Innocuous as it may seem at first sight selecting pairs instead of singletons introduces considerable complexity into the picture and requires a careful analysis which quickly gets intractable. However, in particular cases of simple combinatorial structures (and their generating functions) it is possible to carry all the calculus through. An example that we will consider in detail is the Glaisher-Crofton formula of Eq.~(\ref{Glaisher}) which evaluates action of the exponential of second derivative on the gaussian. Before we proceed to this result, discussed in Sect.~\ref{Proof}, we will illustrate our combinatorial methodology on a simpler case of the action on a monomial and provide a link with combinatorial model of Hermite polynomials.

Let us start with the explicit expression which obtains from expanding the exponential and differentiating the monomial
\begin{eqnarray}\label{Hermite-xn}
\exp \left(\tfrac{1}{2}(yD_x)^2\right)\,x^n&=&\sum_{k=0}^{\lfloor n/2\rfloor}\frac{n!}{2^kk!(n-2k)!}x^{n-2k}y^{2k}\,.
\end{eqnarray}
For $y=i$ it specialises to the Hermite polynomial $He_n(x)$. More generally, we may also write
\begin{eqnarray}\label{Hermite-exp}
\exp \left(\tfrac{1}{2}(yD_x)^2\right)\ \exp\left(xt\right)&=&\exp \left(\tfrac{1}{2}(yt)^2\right)\cdot \exp \left(xt\right)=\exp \left(xt+\tfrac{1}{2}(yt)^2\right)\,,
\end{eqnarray}
which stems from the fact that $\exp(xt)$ is an eigenvector of the derivative operator $D_x$ to eigenvalue $t$. Again, for $y=i$ it is the exponential generating function of Hermite polynomials (cf. Eq.~(\ref{Hermite})). For the purpose at hand we will leave variable $y$ unspecified so to deal with positive integers only (cf. coefficients in Eq.~(\ref{Hermite-xn})).\footnote{It is a typical combinatorial trick to introduce additional labels (or weights) which often allows to get rid of negative or non-integer factors entering multiplicatively in the expressions. Then enumeration of structures proceeds also with respect to this additional label (or weight) which, if needed, can be specified to the required value at the end.} Our aim is to understand these formulas in terms of enumeration of structures.

\begin{figure}[t]
\begin{center}
\includegraphics[width=\textwidth]{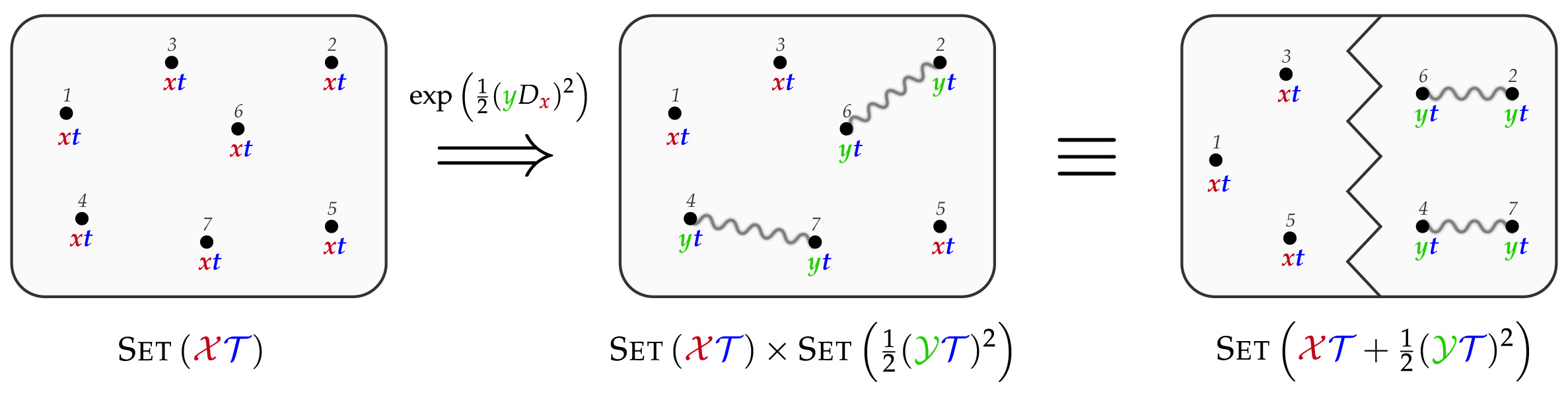}
\end{center}\caption{\label{Fig-Hermite} On the left, an instance of structure in $\mathscr{F}=\set(\mathcal{X}\mathcal{T})$ which is a set of labelled atoms $\mathcal{T}$, each additionally weighted by atom of type $\mathcal{X}$. In the middle, example of a structure arising from the action of $\exp(\tfrac{1}{2}(yD_x)^2)$ on $\mathscr{F}$ which consists in selecting a set of unordered pairs of $\mathcal{X}$ atoms and replacing them by atoms of type $\mathcal{Y}$. Selected pairs are depicted by wavy lines. On the right, decomposition of such structure into a product of two sets comprising singletons $\mathcal{X}\mathcal{T}$ (untouched by derivatives) and doubletons $\tfrac{1}{2}(\mathcal{Y}\mathcal{T})^2$ (selected by derivatives).}
\end{figure}

In order to use combinatorial description of Eq.~(\ref{expyD2}) we interpret $\exp\left(xt\right)$ as the exponential generating function of the labelled class of sets $\set(\mathcal{X}\mathcal{T})$ (see~\cite[Sect. II]{FlSe09} for a precise definition and discussion of labelled classes and their relation with exponential generating functions). It is comprised of sets whose atoms $\mathcal{T}$ carry integer labels $1,2,3,...$, and additionally to each $\mathcal{T}$ is attached an atom (or weight) of type $\mathcal{X}$, see Fig.~\ref{Fig-Hermite} on the left. We have the following translation rule (see Eq.~(\ref{set}))
\begin{eqnarray}
\mathscr{F}=\set(\mathcal{X}\mathcal{T})&\Longrightarrow&F(x,t)=\sum_{n\geqslant0}x^n\tfrac{t^n}{n!}=\exp(xt)\,.
\end{eqnarray} 
(Clearly, for given $n$ one can built one such set and its weight is $x^n$.) Now, following combinatorial description of Eq.~(\ref{expyD2}) action of $\exp \left(\tfrac{1}{2}(yD_x)^2\right)$ on an individual set in $\set(\mathcal{X}\mathcal{T})$ consists in selecting in all possible ways a subset of doubletons. This amounts to splitting each original set into products of two subsets: one comprising singletons and the other doubletons. Additionally, these subsets differ in that atoms in singletons (untouched by derivatives) carry weight $\mathcal{X}$, while each atom forming a doubleton (arising from nontrivial action of second derivatives) carries label $\mathcal{Y}$. Yet another way of seeing the resulting class of objects is to understand them as simply a set of singletons $\mathcal{X}\mathcal{T}$ and doubletons $\tfrac{1}{2}(\mathcal{Y}\mathcal{T})^2$. See Fig.~\ref{Fig-Hermite} for illustration. Formally, one writes the following combinatorial equivalences
\begin{eqnarray}
\exp \left(\tfrac{1}{2}(yD_x)^2\right)\set(\mathcal{X}\mathcal{T})=\set(\mathcal{X}\mathcal{T})\times\set\left(\tfrac{1}{2}(\mathcal{Y}\mathcal{T})^2\right)=\set\left(\mathcal{X}\mathcal{T}+\tfrac{1}{2}(\mathcal{Y}\mathcal{T})^2\right)\,,
\end{eqnarray}
which on the level of generating functions readily transform (cf. Appendix~\ref{Appendix}) into a chain of algebraic equalities providing a combinatorial proof of Eq.~(\ref{Hermite-exp}). As a consequence of this discussion we get a simple combinatorial model of Hermite polynomials. Namely, coefficients of $He_n(x)=\sum_{k=0}^nh_{n,k}x^k$ count the number of possible ways to \emph{select a subset of $k$ doubletons (each weighted by $-1$) out of a set of $n$ distinguishable objects} (it is exactly the coefficient of $x^{n-2k}$ in Eq.~(\ref{Hermite-xn}) for $y=i$).\footnote{Another way to see it directly from our combinatorial description of the exponential of second derivative is to interpret $x^n$ as the generating function of a single $n$ element set, i.e. we have $\mathscr{F}=\set_n(\mathcal{X})\Rightarrow F(x)=x^n$. Then combinatorial model in terms of the choice of doubletons is a simple consequence of the specification of Eq.~(\ref{expyD2}).} We note that this model is a  rephrasing of the interpretation of Hermite polynomials in terms of weighted involutions~\cite[Sect. 2.3.]{BeLaLe98}. Moreover, it can be straightforwardly extended to provide a combinatorial interpretation of a larger class of multivariate Hermite-Kamp\'e de F\'eriet polynomials~\cite{DaOtToVa97,RiTa09}.

\section{Proof of the Glaisher--Crofton identity}
\label{Proof}

Here we prove identity Eq.~(\ref{Glaisher}) by a purely combinatorial argument by analysing structures generated by the second derivative discussed above.
We will proceed in a step by step manner explaining details of combinatorial constructions and structures that appear along the way. Although most of them are standard in combinatorial community we take a rather explicit and methodological route that may be of help for an unaccustomed reader. In the following we adopt standard notation from the book~\cite{FlSe09} (see also Appendix~\ref{Appendix}).

\begin{figure}[t]
\begin{center}
\includegraphics[width=\textwidth]{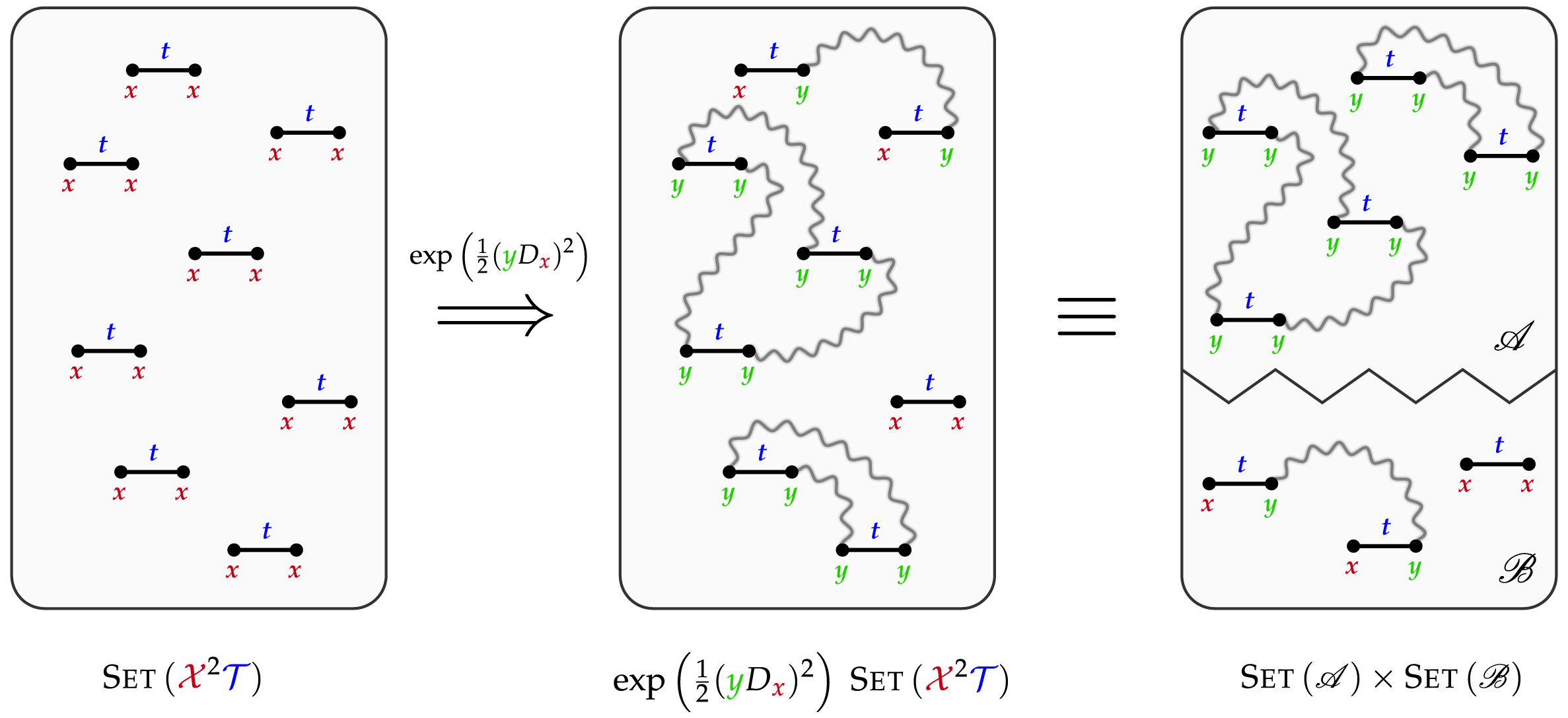}
\end{center}\caption{\label{Fig-Glaisher} On the left, an instance of structure in $\mathscr{C}=\set(\mathcal{X}^2\mathcal{T})$ which is a set of labelled atoms $\mathcal{T}$, each having a finer structure comprised of two atoms of type $\mathcal{X}$ (say the left and the right one). In the middle, example of a structure arising from the action of $\exp(\tfrac{1}{2}(yD_x)^2)$ on $\mathscr{C}$ which consists in selecting a set of unordered pairs of $\mathcal{X}$ atoms and replacing them by atoms of type $\mathcal{Y}$. Selected pairs are depicted by wavy lines. On the right, decomposition of such a structure into a product of two sets comprised of closed ($\mathscr{A}$) and open ($\mathscr{B}$) chains respectively.}
\end{figure}

Our goal is to calculate explicit form of the expression (cf. the left hand side of Eq.~(\ref{Glaisher}) with $\alpha=y^2/2$ and $t=-1$):
\begin{eqnarray}\label{Glaisher-proof}
\exp\left(\tfrac{1}{2}(yD_x)^2\right)\,\exp(x^2t)=\sum_{k\geqslant0}\tfrac{1}{k!}\left(yD_x\right)^{2k}\cdot\sum_{n\geqslant0}\frac{(x^2t)^n}{n!}=R(x,y,t)\,.
\end{eqnarray}
Following our combinatorial strategy we will treat $\exp(x^2t)$ as the exponential generating function enumerating some combinatorial objects. Let us define them as a labelled class of sets $\set(\mathcal{X}^2\mathcal{T})$. This class is comprised of sets of labelled atoms $\mathcal{T}$ (i.e. each atom carries integer label) which are weighted with two atoms of type $\mathcal{X}$. More pictorially, we will depict them as sets of doubletons of unlabelled atoms $\mathcal{X}$ such that each doubleton carries a labelled marker $\mathcal{T}$. See Fig.~\ref{Fig-Glaisher} on the left for illustration. Clearly, we have the following translation rule (cf. Eq.~(\ref{set}))
\begin{eqnarray}
\mathscr{C}=\set(\mathcal{X}^2\mathcal{T})&\Longrightarrow&C(x,t)=\sum_{n\geqslant0}x^{2n}\tfrac{t^n}{n!}=\exp(x^2t)\,.
\end{eqnarray}
We know from discussion of Eq.~(\ref{expyD2}) that on combinatorial level exponential $\exp(\tfrac{1}{2}(yD_x)^2)$ consists in selecting in all possible ways unordered pairs of $\mathcal{X}$'s (not necessarily attached to the same $\mathcal{T}$) and replacing each $\mathcal{X}$ in chosen pairs by $\mathcal{Y}$. Fig.~\ref{Fig-Glaisher} in the middle provides a generic example a structure arising in this procedure. If we denote the resulting class of structures by $\mathscr{R}$, then we may write
\begin{eqnarray}
\mathscr{R}=\exp(\tfrac{1}{2}(yD_x)^2)\,\mathscr{C}&\Longrightarrow&\exp\left(\tfrac{1}{2}(yD_x)^2\right)\,\exp(x^2t)=R(x,y,t)\,,
\end{eqnarray}
where $R(x,y,t)=\sum_{k,l,n\geqslant0}R_{k,l,n}\,x^ky^l\frac{t^n}{n!}$ is the exponential generating function enumerating structures in $\mathscr{R}$. Note that this gives a precise combinatorial meaning to Eq.~(\ref{Glaisher-proof}).

Since we have reduced our goal to finding the exponential generating function $R(x,y,z)$ we need to come up with a systematic specification of structures in $\mathscr{R}$. For this purpose let us observe that doubletons, initially detached one from another in $\mathscr{C}$, now tie together to form either open or closed chains. One can group them together to split each structure in $\mathscr{R}$ into a product of two sets: one containing only closed and the second open chains. This entails the following combinatorial decomposition and its translation to the exponential generating functions (cf. Eqs.~(\ref{product}) and (\ref{set}))
\begin{eqnarray}\label{R}
\mathscr{R}=\set(\mathscr{A})\times\set(\mathscr{B})&\Longrightarrow&R(x,y,t)=\exp (A(x,y,t))\cdot \exp(B(x,y,t))\,,
\end{eqnarray}
where $\mathscr{A}$ is a class of open chains and $\mathscr{B}$ is a class of closed chains, see Fig.~\ref{Fig-Glaisher} on the right. Hence the problem comes down to finding two exponential generating functions $A(x,y,t)$ and $B(x,y,t)$ enumerating objects of type $\mathscr{A}$ and $\mathscr{B}$.

\begin{figure}[t]
\begin{center}
\includegraphics[width=\textwidth]{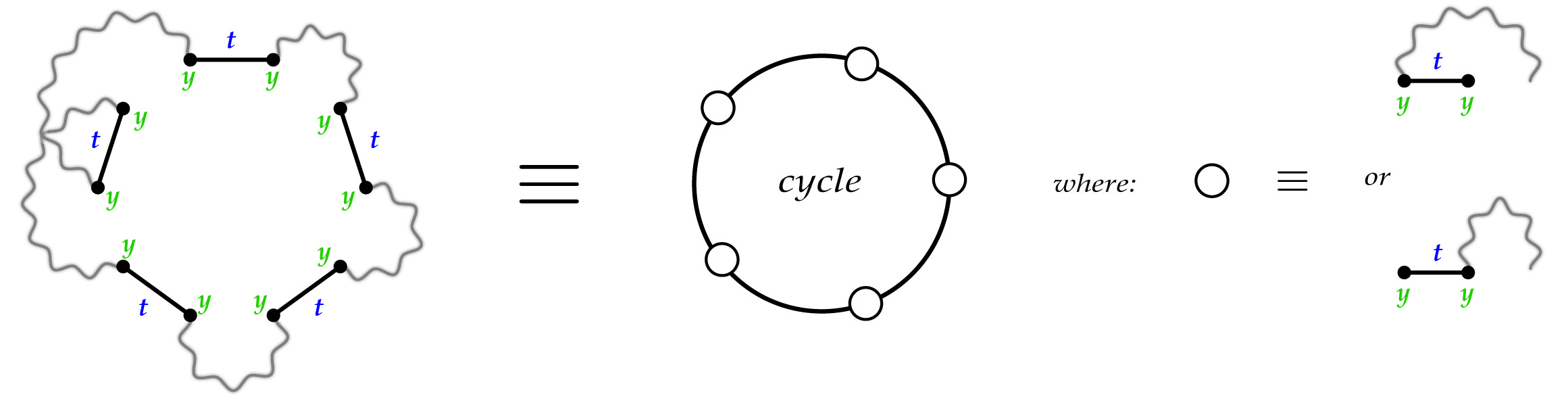}
\end{center}\caption{\label{Fig-Glaisher-A} A generic closed chain in $\mathscr{A}$ ambled in the plane can be seen as a cycle built out of two kinds of blocks, each of type $\mathcal{Y}^2\mathcal{T}$, arising from different choices of $\mathcal{X}$'s attached to the same $\mathcal{T}$ (either left or right) by the derivative (hence possible crossings of the wavy lines).}
\end{figure}

\begin{figure}[t]
\begin{center}
\includegraphics[width=\textwidth]{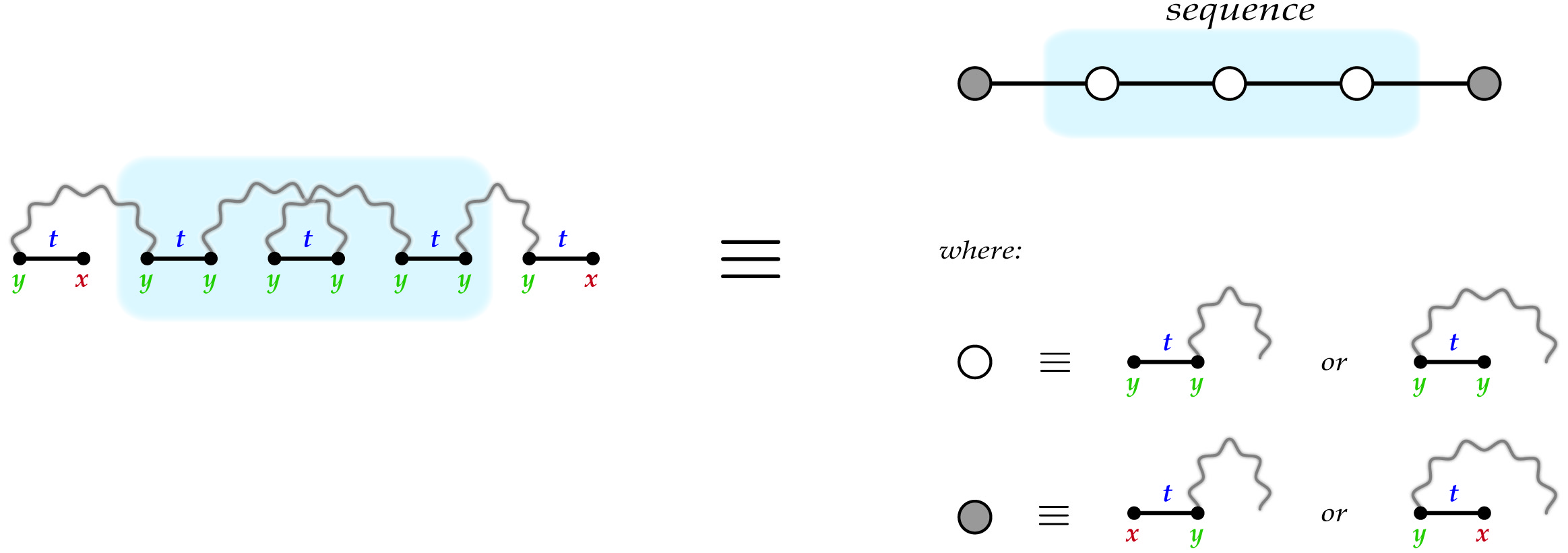}
\end{center}\caption{\label{Fig-Glaisher-B} If embedded in the line, each open chain in $\mathscr{B}$ forms a sequence whose inner part (in the grey box) is built out of blocks of type $\mathcal{Y}^2\mathcal{T}$ with two blocks of type $\mathcal{X}\mathcal{Y}\mathcal{T}$ attached at the ends. Here again, each block can have two configurations due to different choices made by derivatives (possible crossings of the wavy lines).}
\end{figure}

Let us start with the class of closed chains $\mathscr{A}$. Its elements embedded in the plane can be simply seen as cycles whose building blocks have a finer structure of type $\mathcal{Y}^2\mathcal{T}$ which occurs in two possible configurations arising from two possible choices of $\mathcal{X}$ by the derivative in the same doubleton, see Fig.~\ref{Fig-Glaisher-A} for pictorial explanation. Therefore the whole class specifies as follows
\begin{eqnarray}\label{A}
\mathscr{A}=\cyc(2\mathcal{Y}^2\mathcal{T})\,,
\end{eqnarray}
and by the standard transfer rule for labelled classes, Eq.~(\ref{cycle}), we get
\begin{eqnarray}\label{A-gen}
A(x,y,t)=\log\frac{1}{1-2y^2t}\,.
\end{eqnarray}

The class of open chains $\mathscr{B}$ can be described in a similar manner. First we observe that each such chain can be embedded in the line in two possible ways. Then it forms a sequence which can be decomposed into the inner part which is a sequence of blocks of type $\mathcal{Y}\mathcal{T}\mathcal{Y}$ with two additional blocks of type $\mathcal{X}\mathcal{Y}\mathcal{T}$ and $\mathcal{Y}\mathcal{T}\mathcal{X}$ attached at the ends (left and right respectively). Here as well, each block occurs in two possible configurations arising from two possible choices of $\mathcal{X}$ by the derivative in the same doubleton. See Fig.~\ref{Fig-Glaisher-B} for illustration. This gives the following combinatorial specification
\begin{eqnarray}\label{B}
\mathscr{B}=\mathcal{X}\mathcal{T}\mathcal{X}+\tfrac{1}{2}\,(2\,\mathcal{X}\mathcal{Y}\mathcal{T})\times \seq((2\,\mathcal{Y}\mathcal{T}\mathcal{Y}) \times(2\,\mathcal{Y}\mathcal{T}\mathcal{X})\,.
\end{eqnarray}
Note that coefficient $\tfrac{1}{2}$ stems from double counting due to embedding in a line and the term $\mathcal{X}\mathcal{T}\mathcal{X}$ makes up for a single structure left out by the above description. Having specified $\mathscr{B}$ we obtain exponential generating function by means of standard transfer rules, Eq.~(\ref{sequence}), which yield
\begin{eqnarray}\label{B-gen}
B(x,y,t)=x^2t+\tfrac{1}{2}\ 2xyt\cdot\frac{1}{1-2y^2t}\cdot2xyt=\frac{x^2t}{1-2y^2t}\,.
\end{eqnarray}

Now, by substituting Eqs.~(\ref{A-gen}) and (\ref{B-gen}) to Eqs.~(\ref{R}), we get
\begin{eqnarray}
R(x,y,t)=\frac{1}{1-2y^2t}\exp\left(\frac{x^2t}{1-2y^2t}\right)\,.
\end{eqnarray}
This completes derivation of the Glaisher-Crofton identity, Eq.~(\ref{Glaisher}), which readily obtains from Eq.~(\ref{Glaisher-proof}). In conclusion, let us remark that we benefit from the proof by combinatorial insight into the nature of both factors on the r.h.s. of Eq.~(\ref{Glaisher}) which can be interpreted as exponential generating functions enumerating sets of closed and open chains respectively formed by derivatives acting on a gaussian.

\section{Discussion and outlook}

Many computational problems require keeping track and skilful rearrangement of terms involved in algebraic expressions. It often comes down to the analysis of their structural properties and counting terms grouped with respect to some characteristics. This is a natural domain of application for modern combinatorics which has developed a large array of tools for systematic treatment of such problems. In this paper we have considered a few examples where it can be effectively used for evaluation of the action of exponential in derivatives on a function. Fundamental in this approach is treatment of a function as a generating function enumerating some simple combinatorial objects. This shift in perspective allows to take derivatives (and their exponentials) as combinatorial constructors which produce a new class of objects which often can be enumerated with combinatorial flair. We have illustrated this approach by showing simple combinatorial proofs and interpretation of Taylor's formula, connecting exponential in second derivative with a model of Hermite polynomials, and deriving the Glaisher-Crofton identity. It is worth emphasising that the proofs are purely combinatorial, thus do not require any arguments involving integral representations or analyticity. 

Crucial to our development was combinatorial understanding of the exponential of second derivative which consists in selecting a collection of unordered pairs in a structure it acts on. It can be also seen as superposing structure of doubletons on the original one which is connected with the Hadamard product of generating functions considered in various combinatorial contexts (see e.g.~\cite{BeLaLe98,FlSe09,BlPeSoHoDu05}).
We also mention a natural link with a rich framework of umbral calculus~\cite{Ro84} where polynomial sequences can be considered as generated by the action of differential operators (see also monomiality principle~\cite{BlDaHoPe06}) . Such description is attainable for a large family of Sheffer-type polynomials (including binomial-type and Appell sequences) and therefore admit combinatorial interpretation along the lines considered in the present paper. This theme will be the subject of subsequent publication~\cite{TBP-Bl-etal}.

We observe that our discussion is not limited to first and second derivative only. It can be straightforwardly extended to derivatives of higher order which correspond to selecting subsets of higher cardinality. Moreover, one can generalise this framework to partial derivatives in several variables and multivariate polynomials (e.g. Hermite or Kamp\'e-de-F\'eriet polynomials~\cite{DaOtToVa97,RiTa09}) by interpreting them as enumerating combinatorial structures built of atoms of several kinds~\cite{TBP-Bl-etal}. 

Finally, let us remark that combinatorial approach to derivatives also provides interesting insights into operator identities. One example is a systematic treatment of the normal ordering problem~\cite{BlFl11,BlHoPeSoDu07}. Clearly, majority of operator identities admit combinatorial interpretation as they typically arise from algebraic manipulation of discrete structures~\cite{BeLaLe98,LaLa09}. As such it opens the whole field of application for combinatorial approach. In this paper we have illustrated this point only on a few examples which can be seen as instances of a broad class of operator identities amenable to combinatorial methodology (cf. Sack identity, Baker-Campbell-Hausdorff formula, Rodrigues-type formulas, Crofton identities, etc.~\cite{DaOtToVa97,RiTa09,TBP-Bl-etal}).

\section*{Appendix: Combinatorial constructions}\label{Appendix}

Our primary reference for combinatorial analysis is a standard book \emph{Analytic Combinatorics} by Ph. Flajolet and R. Sedgewick~\cite{FlSe09}. 
Here, we briefly recall basic terminology and a few standard translation rules for labelled constructions used in this paper. 

Suppose we are given a \emph{combinatorial class} $\mathscr{C}$ which consists of a denumerable collection of objects built of \emph{labelled} atoms $\mathcal{T}$ (see~\cite[Ch.~II]{FlSe09} for precise definition of labelled class). Usually size of an object is the number of atoms it is built of, and a typical problem is to count the number of structures of a given \emph{size}. In other words, one seeks the sequence $c_n=\#\,\mathscr{C}_n$, where $\mathscr{C}_n=\left\{\phi\in\mathscr{C}:|\phi|=n\right\}$ which is conveniently encoded in the \emph{exponential generating function} (e.g.f.)
\begin{eqnarray}
C(t)=\sum_{n\geqslant0}c_n\,\frac{t^n}{n!}=\sum_{\phi\in\mathscr{C}}\frac{t^{|\phi|}}{|\phi|!}\,.
\end{eqnarray}
Let us remark that the reason for the use of exponential generating functions, rather that ordinary generating functions  (o.g.f.), is simplicity of transfer rules in the domain of labeled classes. (Ordinary generating functions are typically used for enumeration of unlabelled structures, cf.~\cite[Ch.~I]{FlSe09}.)

The point of combinatorial analysis of structures is construction of complex classes from simpler ones. The initial building blocks include the \emph{atomic class} $\mathcal{T}$, which comprises a single element of size $1$ and has e.g.f. $t$, and the \emph{neutral class} $\mathcal{E}$, which consists of a single element of size $0$ and has e.g.f. 1. Then complex structures are built by well defined set theoretical constructions which provide a precise specification of the class. Remarkably, these constructions can be translated into algebraic equations for the corresponding generating functions which solve the enumeration problem. Below, we give a short list of such constructions and translation rules that we exploit in this paper.

The most basic one is \emph{disjoint union}, henceforth denoted by "$+$", which clearly corresponds to
\begin{eqnarray}\label{sum}
\mathscr{C}=\mathscr{A}+\mathscr{B}&\Longrightarrow&C(t)=A(t)+B(t)\,.
\end{eqnarray}
Another one is the \emph{labelled product}, denoted by "$\star$", which forms a cartesian product of objects and relabels the atoms in order-consistent manner. We have the following translation rule
\begin{eqnarray}\label{product}
\mathscr{C}=\mathscr{A}\star\mathscr{B}&\Longrightarrow&C(t)=A(t)\cdot B(t)\,.
\end{eqnarray}
If objects of one structure $\mathscr{B}$ are \emph{substituted} into atoms of another structure $\mathscr{A}$ and relabelled in the order-consistent way, then the e.g.f. of such constructed class is given by (assuming $\mathscr{B}_0=\emptyset$, i.e. $B(0)=0$)
\begin{eqnarray}\label{substitution}
\mathscr{C}=\mathscr{A}\circ\mathscr{B}&\Longrightarrow&C(t)=A\left(B(t)\right)\,.
\end{eqnarray}
It is then possible to form the class of all (labelled) \emph{sequences}, \emph{sets} and \emph{cycles} (respectively denoted by $\set$, $\seq$ and $\cyc$) with components in $\mathscr{A}$. The corresponding generating functions are given by the following dictionary (assuming $\mathscr{A}_0=\emptyset$, i.e. $A(0)=0$)
\begin{eqnarray}
\label{set}
\mathscr{C}=\set(\mathscr{A})&\Longrightarrow&C(t)=\exp(A(t))\,,\\
\label{cycle}
\mathscr{C}=\cyc(\mathscr{A})&\Longrightarrow&C(t)=\log\frac{1}{1-A(t)}\,,\\
\label{sequence}
\mathscr{C}=\seq(\mathscr{A})&\Longrightarrow&C(t)=\frac{1}{1-A(t)}\,.
\end{eqnarray}
The is a non-exhaustive selection of possible constructions which serves the purposes of the present paper. For a comprehensive survey of the methods of combinatorial enumeration via generating functions we refer to the classic books on this subject~\cite{FlSe09,BeLaLe98,Wi06}.


\paragraph{Acknowledgments}

This paper is dedicated to the memory of the late Philippe Flajolet and Allan I. Solomon. We had a great privilege benefit from their exceptional joy and understanding of science.

\bibliographystyle{alpha}

\bibliography{/Users/Blasiak/Research/Library/Bibliography/CombQuant}

\vspace{0.5cm}

\newpage

\noindent{\bf\textsc{Pawel Blasiak \& Andrzej Horzela}}

\noindent H. Niewodnicza\'nski Institute of Nuclear Physics, Polish Academy of Sciences

\noindent ul. Eliasza-Radzikowskiego 152, 31-342 Krak\'ow, Poland

\vspace{0.5cm}

\noindent{\bf\textsc{Gerard H. E. Duchamp}}

\noindent Laboratoire d'Informatique de Paris-Nord, Universit\'e Paris XIII

\noindent Institut Galil\'e, 99 Av. J.-B. Clement, 93430 Villetaneuse, France

\vspace{0.5cm}

\noindent{\bf\textsc{Karol A. Penson}}

\noindent Laboratoire de Physique Th\'eorique de la Mati\'ere Condens\'ee, Universit\'e Paris VI

\noindent Bo\^ite 121, 4 Pl. Jussieu, 75252 Paris Cedex 05, France

\end{document}